\documentclass[12pt]{amsart}
\usepackage{amsmath,amssymb,palatino,cite}

\newcommand{\eps}{\varepsilon}
\newcommand{\R}{\mathbb{R}}

\newcommand{\RN}{{\mathbb{R}^N}}
\newcommand{\RT}{{\mathbb{R}^3}}

\newcommand{\de}{\partial}

\DeclareMathOperator{\cat}{cat}

\renewcommand{\le}{\leqslant}
\renewcommand{\ge}{\geqslant}

\renewcommand{\b }{\beta }

\renewcommand{\d }{\delta }
\newcommand{\vfi}{\varphi}

\newcommand{\n }{\nabla }

\newcommand{\G}{\Gamma}

\newcommand{\A}{{\mathcal A}}

\renewcommand{\H}{{\mathcal H}}

\renewcommand{\P}{{\mathcal P}}

\newcommand{\tf}{\tilde{\phi}}

\numberwithin{equation}{section}
\newtheorem{theorem}{Theorem}[section]
\newtheorem{lemma}[theorem]{Lemma}

\newtheorem{remark}[theorem]{Remark}
\newtheorem{corollary}[theorem]{Corollary}
\renewenvironment{proof}{\noindent{\textbf{Proof\quad}}}{$\hfill\square$\vspace{0.2 cm}\\}
\newenvironment{proofN}{\noindent{\textbf{Proof of Theorem \ref{thN}\quad}}}{$\hfill\square$\vspace{0.2 cm}\\}

\title[Coupled nonlinear Schr\"odinger systems with potentials]
{Coupled nonlinear Schr\"odinger systems with potentials}

\author[Alessio Pomponio]{Alessio Pomponio$^*$}

\address{Dipartimento di Matematica
\newline\indent
Politecnico di Bari
\newline\indent
Via Amendola 126/B, I-70126 Bari, Italy}
\email{a.pomponio@poliba.it}

\thanks{$*$ The author was partially supported by the MIUR research
project ``Metodi Variazionali ed Equazioni Differenziali Nonlineari''}

\begin{document}

\begin{abstract}
Coupled nonlinear Schr\"odinger systems describe some physical phenomena 
such as the propagation in birefringent
optical fibers, Kerr-like photorefractive media in optics and Bose-Einstein
condensates. In this paper, we study the existence of concentrating solutions of 
a singularly perturbed coupled nonlinear Schr\"odinger system, in presence of potentials. 
We show how the location of the concentration points depends strictly on the potentials.
\end{abstract}

\maketitle

\section{Introduction}

Very recently, different authors focused their attention on coupled nonlinear Schr\"odinger systems
which describe physical phenomena such as the propagation in birefringent
optical fibers, Kerr-like photorefractive media in optics and Bose-Einstein
condensates.

First of all, let us recall that, in the last twenty years, 
motivated by the study of the propagation of pulse in
nonlinear optical fiber, the nonlinear Schr\"odinger equation,
\[
- \Delta u + u=u^3 \qquad {\rm in}\;\RT,
\]
has been faced by many authors. It has been 
proved the existence of the least energy solution (ground state solution), which is radial
with respect to some point, positive and exponentially decaying
with its first derivatives at infinity. Moreover there are also many
papers about the semiclassical states for the nonlinear
Schr\"odinger equation with the presence of potentials
\[
- \eps^2 \Delta u + V(x)u=u^3 \qquad {\rm in}\;\RT,
\]
giving sufficient and necessary conditions to the existence of solutions concentrating in
some points, and recently, in set with non-zero dimension, (see e.g.
\cite{AFM,AMN1,AMN2,AMS,BL,dPF,FW,GNN,K,Li,O,Rab,W,WZ}).

However, by I.P.~Kaminow~\cite{Ka}, we know that single-mode
optical fibers are not really ``single-mode'', but actually
bimodal due to the presence of birefringence. This birefringence
can deeply influence the way in which an optical evolves during
the propagation along the fiber. Indeed, it can occur that the
linear birefringence makes a pulse split in two, while
nonlinear birefringent traps them together against splitting.
C.R.~Menyuk \cite{M1,M2} showed that the evolution of two
orthogonal pulse envelopes in birefringent optical fibers is
governed by the following coupled nonlinear Schr\"odinger system
\begin{equation}\label{eq:M}
\begin{cases}
i \phi_t +\phi_{xx}+|\phi|^2 \phi +\b |\psi|^2 \phi=0,
\\
i \psi_t +\psi_{xx}+|\psi|^2 \psi +\b |\phi|^2 \psi=0,
\end{cases}
\end{equation}
with $\b$ positive constant depending on the anisotropy of the
fibers. System \eqref{eq:M} is also important for industrial applications in
fiber communications systems \cite{HK} and all-optical switching devices
\cite{I}. If one seeks for standing wave solutions of \eqref{eq:M},
namely solutions of the form
\[
\phi(x,t)= e^{i w_1^2 t} u(x) \quad {\rm and} \quad \psi(x,t)= e^{i w_2^2 t} v(x),
\]
then \eqref{eq:M} becomes
\begin{equation}\label{eq:M'}
\begin{cases}
-u_{xx} + u =|u|^2 u + \b |v|^2 u  & {\rm in} \; \R,
\\
-v_{xx} + w^2 v =|v|^2 v + \b |u|^2 v  & {\rm in} \; \R,
\end{cases}
\end{equation}
with $w^2=w_2^2 / w_1^2$. 
Finally we want to recall that \eqref{eq:M'} describes also other physical
phenomena, such as Kerr-like photorefractive media in optics, (cf. \cite{AA,CCMS}).

Problem \eqref{eq:M'}, in a more general situation and also in higher dimension, has been studied by 
R.~Cipolatti \& W.~Zumpichiatti \cite{CZ1,CZ2}. By concentration compactness arguments, 
they prove the existence and the regularity of the ground state solutions $(u,v)\neq (0,0)$. 
Later on, in two very recent papers, T.C.~Lin \& J.~Wei \cite{LW1} and 
L.A.~Maia, E.~Montefusco \& B.~Pellacci \cite{MMP} deal with problem
\eqref{eq:M'}, also in the multidimensional case, and, among other results, they prove the existence of
least energy solutions of the type $(u,v)$, with $u,v>0$. Moreover 
T.C.~Lin \& J.~Wei \cite{LW1} prove that, if
$\b<0$, then the ground state solution for \eqref{eq:M'} does not exist. 
We refer to all these papers and to references therein for more complete 
informations about \eqref{eq:M'}.

Another motivation to the study of coupled Schr\"odinger systems arises from the Hartree-Fock
theory for the double condensate, that is a binary mixture of Bose-Einstein
condensates in two different hyperfine states $|1\rangle$ and $|2\rangle$
(cf. \cite{EGBB}). Indeed these phenomena are governed by the following system:
\begin{equation}\label{eqf}
\left\{
\begin{array}{ll}
-\eps^2 \Delta u + \lambda_1 u = \mu_1 u^3 +\b u v^2  & {\rm in} \; \Omega,
\\
-\eps^2 \Delta v + \lambda_2 v = \mu_2 v^3 +\b u^2 v  & {\rm in} \; \Omega,
\\
u,\, v>0 & {\rm in} \; \Omega,
\\
u=v=0 & {\rm on} \; \de\Omega,
\end{array}
\right.
\end{equation}
where $\Omega$ is a bounded domain of $\RT$. Physically, $u$ and $v$ represent 
the corresponding condensate amplitudes,
$\eps^2= \frac{\hbar^2}{2 m}$, with $\hbar$ the Planck constant and $m$ the atom mass. Moreover
$\mu_j = - (N_j -1)U_{jj}$ and $\b= -N_2 U_{12}$, with $N_j \ge 1$ a fixed number of atoms in
the hyperfine state $|j\rangle$, and $U_{ij}= 4 \pi \frac{\hbar^2}{m} a_{ij}$, where $a_{jj}$'s and
$a_{12}$ are the intraspecies and interspecies scattering lengths. Besides, by
E.~Timmermans \cite{T}, we infer that $\mu_j=\mu_j(x)$ represents a chemical potential.
For more informations about \eqref{eqf}, see \cite{LW1,LW2} and references therein.

T.C.~Lin \& J.~Wei, in \cite{LW2}, studied problem \eqref{eqf} with $\lambda_1,\,\lambda_2,\,\mu_1,\,\mu_2$
positive constant and they proved that if $\b < \sqrt{\mu_1 \mu_2}$, for $\eps$ sufficiently small,
\eqref{eqf} has a least energy solution $(u_\eps,v_\eps)$. Moreover, they distinguished two cases:
the attractive case and the repulsive one. In the attractive case, which occurs whenever $\b>0$,
$u_\eps$ and $v_\eps$ concentrate respectively in $Q_\eps$ and $Q'_\eps$, with
\[
\frac{|Q_\eps -Q'_\eps|}{\eps}\to 0, \qquad {\rm as}\; \eps \to 0.
\]
Precisely they proved that
\begin{align*}
& d(Q_\eps, \de \Omega) \to \max_{Q \in \Omega} d(Q, \de \Omega),
\\
& d(Q'_\eps, \de \Omega) \to \max_{Q \in \Omega} d(Q, \de \Omega).
\end{align*}
In the repulsive case, that is when $\b<0$, the concentration points $Q_\eps$ and $Q'_\eps$
satisfy the following condition:
\[
\vfi(Q_\eps,Q'_\eps) \to \max_{(Q,Q')\in \Omega^2} \vfi(Q,Q'),
\]
where
\[
\vfi(Q,Q')=
\min \{ \sqrt{\lambda_1}|Q-Q'|,\sqrt{\lambda_2}|Q-Q'| , \sqrt{\lambda_1}d(Q, \de \Omega), \sqrt{\lambda_2}d(Q', \de \Omega)\}.
\]
In particular
\[
\frac{|Q_\eps -Q'_\eps|}{\eps}\to \infty, \qquad {\rm as}\; \eps \to 0.
\]

Motivated by these results and by the fact that we know that $\mu_j$ may be
not constants (cf. \cite{T}), in this paper we study the following problem:
\begin{equation}\label{eq}\tag{${\mathcal P}_\eps$}
\left\{
\begin{array}{ll}
-\eps^2 \Delta u + J_1(x) u = J_2(x) u^3 +\b u v^2  & {\rm in} \; \Omega,
\\
-\eps^2 \Delta v + K_1(x) v = K_2(x) v^3 +\b u^2 v  & {\rm in} \; \Omega,
\\
u,\, v>0 & {\rm in} \; \Omega,
\\
u=v=0 & {\rm on} \; \de\Omega,
\end{array}
\right.
\end{equation}
with $\Omega \subset \RT$, possibly unbounded and with smooth boundary, and with
$\b<0$, namely in the repulsive case. We will show that the presence of the potentials change
drastically the situation with respect to the case with positive constants for
what concerns the location of peaks, but, in some sense, not the repulsive 
nature of the problem. In fact,
with suitable assumptions on the potentials, for $\eps$
sufficiently small, we will find solutions $(u_\eps, v_\eps)$ of \eqref{eq}, even if not of
least energy, concentrating respectively on $Q_\eps$ and $Q_\eps'$ which tend
toward the same point, determined by the potentials, as $\eps \to 0$, 
but with the property that the distance between them divided by $\eps$
diverges (see Remark \ref{re:QQ})

Up to our knowledge, in this paper we give a first existence result of 
concentrating solutions for problem \eqref{eq}, in presence of potentials.

On the potentials $J_i$ and $K_i$ we will do the following assumptions:
\begin{itemize}
\item[{\bf (J)}] for $i=1,2$, $J_i\in C^1 (\Omega, \R)$, $J_i$ and $D J_i$ are bounded; moreover,
\begin{equation*}
J_i(x)\ge C>0 \quad \textrm{for all } x\in\Omega;
\end{equation*}
\item[{\bf (K)}] for $i=1,2$, $K_i\in C^1 (\Omega, \R)$, $K_i$ and $D K_i$ are bounded; moreover,
\begin{equation*}
K_i(x)\ge C>0 \quad \textrm{for all } x\in\Omega.
\end{equation*}
\end{itemize}

Without lost of generality, we can suppose that there exists $\eps_0>0$, such that
$\Omega_0:= \Omega \cap (\Omega - \eps_0 \,e_1)\neq \emptyset $, where $e_1=(1,0,0)$.

Let us introduce an auxiliary function which will play a crucial role in the study of \eqref{eq}.
Let $\G\colon \Omega_0 \to \R$ be a function so defined:
\begin{equation}\label{eq:Gamma}
\G(Q)=J_1(Q)^{\frac 12}J_2(Q)^{-1} + K_1(Q)^{\frac 12}K_2(Q)^{-1}.
\end{equation}
Let us observe that by {\bf (J)} and {\bf (K)}, $\G$ is well defined.

Our main result is:

\begin{theorem}\label{th1}
Suppose {\bf (J)} and {\bf (K)} and $\b <0$. Let  $Q_0 \in \Omega_0$ be an
isolated local strict minimum or maximum of $\G$.
There exists $\bar \eps>0$ such that if $0<\eps<\bar \eps$,
then \eqref{eq} possesses a solution $(u_\eps,v_\eps)$ such that $u_\eps$
concentrates in $Q_\eps$ with $Q_\eps \to Q_0$, as $\eps \to 0$, and $v_\eps$
concentrates in $Q'_\eps$ with $Q'_\eps \to Q_0$, as $\eps \to 0$.
\end{theorem}

\begin{remark}\label{re:QQ}
Let us observe that, by the proof, it will be clear that, even if
\[
|Q_\eps - Q'_\eps| \to 0, \qquad {\rm as }\; \eps \to 0,
\]
we have
\[
\frac{|Q_\eps - Q'_\eps|}{\eps} \to \infty, \qquad {\rm as }\; \eps \to 0.
\]
\end{remark}

Let us present how Theorem \ref{th1} becomes in 
some particular situations.

Let $H\colon \Omega\to \R$ satisfying the assumption:
\begin{itemize}
\item[{\bf (H)}] $H\in C^1 (\Omega, \R)$, $H$ and $D H$ are bounded; moreover,
\begin{equation*}
H(x)\ge C>0 \quad \textrm{for all } x\in\Omega.
\end{equation*}
\end{itemize}

\begin{corollary}
Suppose {\bf (H)} and $\b <0$. Suppose, moreover, that we are in one of 
the following situations:
\begin{itemize}
\item all the potentials $J_i$ and $K_i$ coincide with $H$;
\item there exists $i_0=1,2$ such that $J_{i_0}\equiv H$ and $K_{i_0}\equiv H$, for $i=i_0$, 
while $J_i$ and $K_i$ are constant for  $i\neq i_0$;
\item all the potentials $J_i$ and $K_i$ are constant, except only one, which coincides with $H$.
\end{itemize}
Let  $Q_0 \in \Omega_0$ be an isolated local strict minimum or maximum of $H$.
There exists $\bar \eps>0$ such that if $0<\eps<\bar \eps$,
then \eqref{eq} possesses a solution $(u_\eps,v_\eps)$ such that $u_\eps$
concentrates in $Q_\eps$ with $Q_\eps \to Q_0$, as $\eps \to 0$, and $v_\eps$
concentrates in $Q'_\eps$ with $Q'_\eps \to Q_0$, as $\eps \to 0$.
\end{corollary}

\begin{remark}
If, instead of $\b$ constant, we consider $\b \in C^1(\Omega,\R)$, bounded and
bounded above by a negative constant, then we have exactly the same results.
\end{remark}

Finally, we want to observe that we can treat also a more general problem 
than \eqref{eq}. Let us consider, indeed,
\begin{equation}\label{eqN}\tag{$\bar{\mathcal P}_\eps$}
\left\{
\begin{array}{ll}
-\eps^2 \Delta u + J_1(x) u = J_2(x) u^{2p -1} +\b u^{p -1} v^p  & {\rm in} \; \Omega,
\\
-\eps^2 \Delta v + K_1(x) v = K_2(x) v^{2p -1} +\b u^p v^{p -1}  & {\rm in} \; \Omega,
\\
u,\, v>0 & {\rm in} \; \Omega,
\\
u=v=0 & {\rm on} \; \de\Omega,
\end{array}
\right.
\end{equation}
with $\Omega \subset \RN$, possibly unbounded and with smooth boundary, 
$N\ge 3$, $2<2p<2N/(N-2)$ and with $\b<0$.

Also in this case, without lost of generality, we can suppose that there exists $\eps_0>0$, such that
$\bar\Omega_0:= \Omega \cap (\Omega - \eps_0 \,\bar e_1)\neq \emptyset $, where $\bar e_1=(1,0,\ldots,0)\in \RN$.

Let us define now $\bar\G\colon \bar\Omega_0 \to \R$ be a function so defined:
\[
\bar\G(Q)=J_1(Q)^{\frac{p}{p-1}-\frac N2}J_2(Q)^{-\frac{1}{p-1}} 
+ K_1(Q)^{\frac{p}{p-1}-\frac N2}K_2(Q)^{-\frac{1}{p-1}}. 
\]

In this case, Theorem \ref{th1} becomes:

\begin{theorem}\label{thN}
Let $N\ge 3$ and $2<2p<2N/(N-2)$. 
Suppose {\bf (J)} and {\bf (K)} and $\b <0$. Let  $Q_0 \in \bar\Omega_0$ be an
isolated local strict minimum or maximum of $\bar\G$.
There exists $\bar \eps>0$ such that if $0<\eps<\bar \eps$,
then \eqref{eqN} possesses a solution $(u_\eps,v_\eps)$ such that $u_\eps$
concentrates in $Q_\eps$ with $Q_\eps \to Q_0$, as $\eps \to 0$, and $v_\eps$
concentrates in $Q'_\eps$ with $Q'_\eps \to Q_0$, as $\eps \to 0$.
\end{theorem}

\begin{remark}
Let us observe that, if $p=2$ and $N=3$, then Theorem \ref{th1} is nothing else than a particular 
case of Theorem \ref{thN}. Nevertheless, since problem \eqref{eq} is more natural 
and more important by a physical point of view, we prefer to present Theorem \ref{th1} as 
our main result and to prove it directly, showing how, with slight modifications, 
the proof of Theorem \ref{thN} follows.
\end{remark}

Theorem \ref{th1} will be proved as a particular case of a multiplicity result in Section 5 
(see Theorem \ref{th:cat}). 
The proof of the theorem relies on a finite dimensional 
reduction, precisely on the perturbation technique developed in \cite{AB, ABC, AMS}. 
In Section 2 we give some preliminary lemmas and some estimates which will be 
useful in Section 3 and Section 4, 
where we perform the Liapunov-Schmidt reduction, making also the asymptotic 
expansion of the finite dimensional functional. 
Finally, in Section 5, we give also a short proof of Theorem \ref{thN}.

\

\begin{center}{\bf Notation}\end{center}
\begin{itemize} 
\item We denote $\Omega_0:= \Omega \cap (\Omega - \eps_0 \, e_1)$, where $e_1=(1,0,0)$ 
and $\eps_0$ is sufficiently small such that $\Omega_0 \neq \emptyset$.

\item If $r>0$ and $x_0 \in \RT$, $B_r (x_0):= \left\{ x\in\RT : |x- x_0| <r \right\}$.
We denote with $B_r$ the ball of radius $r$ centered in the origin.

\item If $u \colon \RT \to \R$ and $P\in \RT$, we set $u_P := u(\cdot - P)$.

\item If $\eps>0$, we set $\Omega_\varepsilon := \Omega/\eps =\{x\in \RT : \eps x \in \Omega\}$.

\item We denote $\H_\eps=H^1_0(\Omega_\varepsilon) \times H^1_0(\Omega_\varepsilon)$.

\item If there is no misunderstanding, we denote with $\|\cdot \|$ and with $(\cdot \mid \cdot)$
respectively the norm and the scalar product both of $H^1_0(\Omega_\varepsilon)$ and of $\H_\eps$.
While we denote with $\|\cdot \|_{\RT}$ and with $(\cdot \mid \cdot)_{\RT}$
respectively the norm and the scalar product of $H^1(\RT)$.

\item With $C_i$ and $c_i$, we denote generic positive constants, which may
  also vary from line to line.
\end{itemize}

\section{Some preliminary}

Performing the change of variable $x \mapsto \eps x$, problem \eqref{eq} becomes:
\begin{equation}\label{eqe}
\left\{
\begin{array}{ll}
- \Delta u + J_1(\eps x) u = J_2(\eps x)u^3 +\b u v^2 =0 & {\rm in} \; \Omega_\varepsilon,
\\
- \Delta v + K_1(\eps x) v = K_2(\eps x)v^3 +\b u^2 v =0 & {\rm in} \; \Omega_\varepsilon,
\\
u,\, v>0 & {\rm in} \; \Omega_\varepsilon,
\\
u=v=0 & {\rm on} \; \de\Omega_\varepsilon,
\end{array}
\right.
\end{equation}
where $\Omega_\eps=\eps^{-1}\Omega$. Of course if $(u,v)$ is a solution of \eqref{eqe},
then $(u(\cdot/\eps),v(\cdot/\eps))$ is a solution of \eqref{eq}.

Solutions of \eqref{eqe} will be found in 
\[
\H_\eps=H^1_0(\Omega_\varepsilon) \times H^1_0(\Omega_\varepsilon),
\]
endowed with the following norm:
\[
\|(u,v)\|^2_{\H_\eps}=\|u\|^2_{H^1_0(\Omega_\varepsilon)} + \|v\|^2_{H^1_0(\Omega_\varepsilon)}, 
\quad \textrm{ for all }(u,v)\in \H_\eps.
\]
If there is no misunderstanding, we denote with $\|\cdot \|$ and with $(\cdot \mid \cdot)$
respectively the norm and the scalar product both of $H^1(\Omega_\varepsilon)$ and of $\H_\eps$. 

Solutions of \eqref{eqe} are critical points of the functional
$f_\eps \colon \H_\eps \to \R$, defined as
\begin{align*}
f_\eps (u,v) &=
\frac{1}{2}\int_{\Omega_\varepsilon} |\nabla u|^2 +
\frac{1}{2}\int_{\Omega_\varepsilon} J_1(\eps x)u^2 
-\frac{1}{4}\int_{\Omega_\varepsilon} J_2(\eps x)u^{4} 
\\
&\quad+ \frac{1}{2}\int_{\Omega_\varepsilon} |\nabla v|^2 +
\frac{1}{2}\int_{\Omega_\varepsilon} K_1(\eps x) v^2 
-\frac{1}{4}\int_{\Omega_\varepsilon} K_2(\eps x)v^{4} 
\\
&\quad - \frac{\b}{2} \int_{\Omega_\varepsilon} u^2 v^2 .
\end{align*}
If we define $f_\eps^J \colon H^1_0(\Omega_\varepsilon) \to \R$ and $f_\eps^K \colon H^1_0(\Omega_\varepsilon) \to \R$ as
\begin{align*}
f_\eps^J (u) =&\;
\frac{1}{2}\int_{\Omega_\varepsilon} |\nabla u|^2 +
\frac{1}{2}\int_{\Omega_\varepsilon} J_1(\eps x)u^2 
-\frac{1}{4}\int_{\Omega_\varepsilon} J_2(\eps x)u^{4} ,
\\
f_\eps^K (v) =&\;
\frac{1}{2}\int_{\Omega_\varepsilon} |\nabla v|^2 +
\frac{1}{2}\int_{\Omega_\varepsilon} K_1(\eps x)v^2 
-\frac{1}{4}\int_{\Omega_\varepsilon} K_2(\eps x)v^{4} ,
\end{align*}
we have
\[
f_\eps (u,v)= f_\eps^J (u) + f_\eps^K (v) - \frac{\b}{2} \int_{\Omega_\varepsilon} u^2 v^2 .
\]
Furthermore, for any fixed $Q\in \Omega$, we define the two functionals $F^J_Q\colon H^1(\RT) \to \R$
and $F^K_Q\colon H^1(\RT) \to \R$, as follows
\begin{align*}
F^J_Q (u) =&\;
\frac{1}{2}\int_{\RT} |\nabla u|^2 +
\frac{1}{2}\int_{\RT} J_1(Q) u^2 
-\frac{1}{4}\int_{\RT} J_2(Q) u^{4} ,
\\
F^K_Q (v) =&\;
\frac{1}{2}\int_{\RT} |\nabla v|^2 +
\frac{1}{2}\int_{\RT} K_1(Q) v^2 
-\frac{1}{4}\int_{\RT} K_2(Q) v^{4} .
\end{align*}

The solutions of \eqref{eqe} will be found near $(U^Q, V^Q)$, properly
truncated, where $U^Q$ is the unique solution of
\begin{equation}\label{eq:UQ}
\left\{
\begin{array}{ll}
-\Delta  u+ J_1(Q) u= J_2(Q)u^3 &   \text{in }\RT,
\\
u >0 &   \text{in }\R^{3},
\\
u(0)=\max_{\R^{3}} u,
\end{array}
\right.
\end{equation}
and $V_{Q}$ is the unique solution of
\begin{equation}\label{eq:VQ}
\left\{
\begin{array}{ll}
-\Delta  v+ K_1(Q) v= K_2(Q)v^3 &   \text{in }\RT,
\\
v >0 &   \text{in }\R^{3},
\\
v(0)=\max_{\R^{3}} v,
\end{array}
\right.
\end{equation}
for an appropriate choice of $Q \in \Omega_0$. It is easy to see that
\begin{align}
U^Q(x)&=\sqrt{J_1(Q)/J_2(Q)} \cdot W\left(\sqrt{J_1(Q)} \cdot x   \right),\label{eq:UQ1}
\\
V^{Q}(x)&=\sqrt{K_1(Q)/K_2(Q)} \cdot W\left(\sqrt{K_1(Q)} \cdot x \right),\label{eq:VQ1}
\end{align}
where $W$ is the unique solution of
\begin{equation}\label{eq:W}
\left\{
\begin{array}
[c]{lll}
-\Delta  z+ z= z^3 &   \text{in }\R^{3},
\\
z>0 &   \text{in }\R^{3},
\\
z(0)=\max_{\R^{3}}z,
\end{array}
\right.
\end{equation}
which is radially symmetric and decays exponentially at infinity with its
first derivatives (cf. \cite{GNN,K}).

For all $Q\in \Omega_0$, we set $Q'= Q'(\eps, Q)= Q+\sqrt{\eps}\,e_1 \in \Omega$ 
and moreover we call $P=P(\eps, Q)=Q/\eps \in \Omega_\varepsilon$ and $P'=P'(\eps, Q)=Q'/\eps \in \Omega_\varepsilon$. Let us
observe that
\begin{equation}\label{eq:PP'}
|P-P'|=\frac{1}{\sqrt{\eps}} \to 0, \qquad {\rm as}\; \eps \to 0.
\end{equation}

Let $\chi \colon \RT \to \R$ be a smooth function such that
\begin{equation}\label{eq:chi}
\begin{array}{rcll}
\chi(x)\!\!\!\!& =&\!\!\!\! 1, &\quad \hbox{ for } |x| \le \eps^{-1/4};
\\
\chi(x) \!\!\!\!&=&\!\!\!\! 0, &\quad \hbox{ for } |x| \ge 2 \eps^{-1/4};
\\
0 \le \chi (x) \!\!\!\!&\le&\!\!\!\! 1, &\quad \hbox{ for } \eps^{-1/4} \le |x| \le 2 \eps^{-1/4};
\\
|\nabla\chi(x)|\!\!\!\!&\le&\!\!\!\! 2\eps^{1/4}, &\quad \hbox{ for }\eps^{-1/4}\le |x|\le 2\eps^{-1/4}.
\end{array}
\end{equation}
We denote
\begin{align}
U_P(x) :=&\; \chi (x -P) \, U^Q (x-P), \label{eq:U}
\\
V_{P'}(x) :=&\; \chi (x -P')\, V^{Q} (x-P'). \label{eq:V}
\end{align}
Let us observe that $(U_P,V_{P'}) \in \H_\eps$. For $Q$
varying in $\Omega_0$,  $(U_P,V_{P'})$ describes a $3$-dimensional manifold, namely,
\begin{equation}\label{eq:Z}
Z^\eps = \left\{(U_P,V_{P'}): Q \in \Omega_0  \right\}.
\end{equation}

\begin{remark}
Of course, if $\Omega=\RT$, then $\Omega_0=\RT$ and we do not need to truncate $U^Q$ and $V^Q$. 
In this case, we would have simply $U_P=U^Q(\cdot -P)$ and $V_{P'}=V^Q(\cdot -P')$.
\end{remark}

First of all let us give the following estimate which will be very useful in the sequel.
\begin{lemma}\label{le:UV}
For all $Q \in \Omega_0$ and for all $\eps$ sufficiently small, 
if $Q'=Q+\sqrt{\eps}\,e_1$, $P=Q/\eps \in \Omega_\varepsilon$ and
$P'=Q'/\eps \in \Omega_\varepsilon$, then
\begin{equation}\label{eq:UV}
\int_{\Omega_\varepsilon} U^2_P V^2_{P'}  =o(\eps).
\end{equation}
\end{lemma}

\begin{proof}
Let us start observing that, since
\[
|P-P'|=\eps^{-1/2} > 4\eps^{-1/4},
\]
we infer that
\[
B_{2\eps^{-1/4}}(P)\cup B_{2\eps^{-1/4}}(P')=\emptyset.
\]
Therefore, by
the definitions of \eqref{eq:U} and \eqref{eq:V} and by the exponential decay 
of $U_P$ and $V_{P'}$, we get
\begin{align*}
\int_{\Omega_\varepsilon} U^2_P V^2_{P'}   \le & \;
\int_{B_{2\eps^{-1/4}}(P)\cup B_{2\eps^{-1/4}}(P')}\!\!\!\big(U^Q \big)^2 (x-P) \big(V^{Q} \big)^2 (x-P')
\\
\le& \; c_1 \int_{\RT \setminus B_{2\eps^{-1/4}}(P')}\big(V^{Q} \big)^2 (x-P') 
\\
&+ c_2 \int_{\RT \setminus B_{2\eps^{-1/4}}(P)}\big(U^{Q} \big)^2 (x-P) 
=o(\eps).
\end{align*}
This concludes the proof.
\end{proof}

In the next lemma we show that the $3$-dimensional manifold $Z_\eps$, defined
in \eqref{eq:Z}, is actually a manifold of almost critical points of $f_\eps$.

\begin{lemma}\label{le:nf}
For all $Q \in \Omega_0$ and for all $\eps$ sufficiently small, if 
$Q'=Q+\sqrt{\eps}\,e_1$, $P=Q/\eps \in \Omega_\varepsilon$ and
$P'=Q'/\eps \in \Omega_\varepsilon$, then
\begin{equation}\label{eq:nf}
\|\nabla f_\eps(U_P, V_{P'})\|=O(\eps^{1/2}).
\end{equation}
\end{lemma}

\begin{proof}
For all $(u,v)\in \H_\eps$, we have:
\begin{align}
(\nabla f_\eps(U_P,V_{P'}) \mid (u,v)) =&
\int_{\Omega_\varepsilon} \left[
\n U_P \cdot \n u
+ J_1(\eps x) U_P u
- J_2(\eps x) U_P^3 u
\right]  \nonumber
\\
&+\int_{\Omega_\varepsilon} \left[
\n V_{P'} \cdot \n v
+ K_1(\eps x) V_{P'} v
-K_2(\eps x) V_{P'}^3 v
\right]  \nonumber
\\
&-\b \int_{\Omega_\varepsilon} U_P V_{P'}^2 u
-\b \int_{\Omega_\varepsilon} U_P^2 V_{P'} v .  \label{eq:nf-tot}
\end{align}
Let us study the first integral of the right hand side of \eqref{eq:nf-tot}. 
By the exponential decay of $U^Q$ and recalling
that $U^Q$ is solution of \eqref{eq:UQ}, we get
\begin{align}
\int_{\Omega_\varepsilon}  \big[
\n U_P & \cdot \n u 
+ J_1(\eps x) U_P u
- J_2(\eps x)U_P^3 u
\big]  \nonumber
\\
=&
\int_{(\Omega -Q)/\eps \,\cap B_{\eps^{-1/4}}} \left[
\n U^Q \cdot \n u_{-P}
+ J_1(\eps x+Q) U^Q u_{-P} \right]  \nonumber
\\
& - \int_{(\Omega -Q)/\eps \,\cap B_{\eps^{-1/4}}} J_2(\eps x+Q)\left(U^Q\right)^3 u_{-P}\, 
+o(\eps)\nonumber
\\
=&
\int_{\RT} \left[
\n U^Q \cdot \n u_{-P}
+ J_1(\eps x+Q) U^Q u_{-P}
- J_2(\eps x+Q)\left(U^Q\right)^3 u_{-P}
\right]  +o(\eps)\nonumber
\\
=&
\int_{\RT} \left[
\n U^Q \cdot \n u_{-P}
+ J_1(Q) U^Q u_{-P}
- J_2(Q) \left(U^Q\right)^3 u_{-P}
\right]    \nonumber
\\
&+\int_{\RT} \big( J_1(\eps x+Q)- J_1(Q) \big) U^Q u_{-P} \nonumber
\\
&-\int_{\RT} \big( J_2(\eps x+Q)- J_2(Q) \big) \left(U^Q\right)^3 u_{-P}
+o(\eps)\nonumber
\\
=&
\int_{\RT} \big(J_1(\eps x+Q)- J_1(Q) \big) U^Q u_{-P}\nonumber
\\
&-\int_{\RT} \big( J_2(\eps x+Q)- J_2(Q) \big) \left(U^Q\right)^3 u_{-P}
+o(\eps). \label{eq:nf-U}
\end{align}
Moreover, from the assumption $D J_i$ bounded, we infer that
\[
|J_i(\eps x +Q)-J_i(Q)| \le c_1 \eps |x|,
\]
and so,
\begin{align}
\int_{\RT} \big(J_1(\eps x +Q)-J_1(Q) \big) U^Q u_{-P}
\le& \;\|u\| \left( \int_{\RT} |J_1(\eps x +Q)-J_1(Q)|^2 |U^Q|^2 \right)^{1/2}
\nonumber
\\
\le& \;c_1 \|u\| \left(\int_\RT \!\! \eps^2 |x|^2|U^Q|^2  \right)^{1/2}
=O(\eps)\|u\|.  \label{eq:restoJ1}
\end{align}
Analogously,
\begin{equation}\label{eq:restoJ2}
\int_{\RT} \big( J_2(\eps x+Q)- J_2(Q) \big) \left(U^Q\right)^3 u_{-P}
=O(\eps)\|u\|.
\end{equation}
Therefore, by \eqref{eq:nf-U}, \eqref{eq:restoJ1} and \eqref{eq:restoJ2}, we infer
\begin{equation}\label{eq:nf-Ufin}
\int_{\Omega_\varepsilon} \left[
\n U_P \cdot \n u
+ J_1(\eps x) U_P u
- J_2(\eps x) U_P^3 u
\right]  =O(\eps)\|u\|.
\end{equation}
Similarly, since $V^Q$ is solution of \eqref{eq:VQ}, we get
\begin{align}
\int_{\Omega_\varepsilon} \big[
\n V_{P'} \cdot \n v &
+ K_1(\eps x) V_{P'} v
- K_2(\eps x)V_{P'}^3 v
\big]  \nonumber
\\
=&
\int_{\RT} \big(K_1(\eps x+Q+\sqrt{\eps} \,e_1)- K_1(Q) \big) V^Q v_{-P'} \nonumber
\\
&-\int_{\RT} \big(K_2(\eps x+Q+\sqrt{\eps}\, e_1)- K_2(Q) \big) \left(V^Q\right)^3 v_{-P'}
+o(\eps). \label{eq:nf-V}
\end{align}
Therefore, from the assumption $D K_i$ bounded, we infer that
\[
|K_i(\eps x +Q+\sqrt{\eps}\,e_1)-K_i(Q)| \le c_2 \sqrt{\eps} \,|\sqrt{\eps}\, x+e_1 |,
\]
and so,
\begin{align}
\int_{\RT} \big(K_1(\eps x +Q+ &\sqrt{\eps}\, e_1)-K_1(Q) \big) V^Q v_{-P'} \nonumber
\\
&\le \|v\| \left( \int_{\RT} |K_1(\eps x +Q+\sqrt{\eps}\, e_1)-K_1(Q)|^2 |V^Q|^2 \right)^{1/2}
\nonumber
\\
&\le  c_2 \|v\| \left(\int_\RT \eps \,|\sqrt{\eps} \,x+e_1|^2|V^Q|^2  \right)^{1/2} =
O(\eps^{1/2}) \|v\|.  \label{eq:restoK1}
\end{align}
Analogously,
\begin{equation}\label{eq:restoK2}
\int_{\RT} \big( K_2(\eps x+Q+\sqrt{\eps}\, e_1)- K_2(Q) \big) \left(V^Q\right)^3 v_{-P'}
=O(\eps^{1/2})\|v\|.
\end{equation}
Therefore, by \eqref{eq:nf-V}, \eqref{eq:restoK1} and \eqref{eq:restoK2}, we infer
\begin{equation}\label{eq:nf-Vfin}
\int_{\Omega_\varepsilon} \left[
\n V_{P'} \cdot \n v
+ K_1(\eps x) V_{P'} v
- K_2(\eps x) V_{P'}^3 v
\right]  =O(\eps^{1/2})\|v\|.
\end{equation}
Let us study the last two terms of \eqref{eq:nf-tot}. Arguing as in Lemma \ref{le:UV}, we get
\begin{equation}\label{eq:nf-UV2}
\left| \int_{\Omega_\varepsilon} U_P V_{P'}^2 u  \right|
\le c_3\left( \int_{\Omega_\varepsilon} U_P^{4/3} V_{P'}^{8/3} \right)^{3/4} \|u\|
=o(\eps)\|u\|,
\end{equation}
and
\begin{equation}\label{eq:nf-U2V}
\left| \int_{\Omega_\varepsilon} U_P^2 V_{P'} v \right|
=o(\eps)\|v\|.
\end{equation}
Now the conclusion of the proof easily follows by \eqref{eq:nf-tot}, 
\eqref{eq:nf-Ufin}, \eqref{eq:nf-Vfin}, \eqref{eq:nf-UV2} and \eqref{eq:nf-U2V}.
\end{proof}

\section{Invertibility of $D^2 f_\eps$ on $\left(T_{(U_P,V_{P'})}Z^\eps \right)^\perp$}

In this section we will show that $D^2 f_\eps$ is invertible on
$\left(T_{(U_P,V_{P'})}Z^\eps \right)^\perp$, where $T_{(U_P,V_{P'})} Z^\eps$
denotes the tangent space to $Z^\eps$ at the point $(U_P,V_{P'})$.

Let $L_{\eps,Q}:(T_{(U_P,V_{P'})}Z^\eps)^\perp\to
(T_{(U_P,V_{P'})}Z^\eps)^\perp$ denote the operator defined by setting
$(L_{\eps,Q} (h,h') \mid (k,k'))= D^2 f_\eps(U_P,V_{P'})[(h,h'),(k,k')]$.

\begin{lemma}\label{le:inv}
Given $\mu>0$, there exists $C>0$ such that, for $\eps$ small enough 
and for all $Q\in \Omega_0$ with $|Q|\le \mu$, one has that
\begin{equation}\label{eq:inv}
\|L_{\eps,Q} (h,h') \|\ge C \|(h,h')\|,\qquad \forall\; (h,h')\in(T_{(U_P,V_{P'})}Z^{\eps})^{\perp}.
\end{equation}
\end{lemma}

\begin{proof}
First of all, let us observe that, for all $(h,h'),(k,k') \in \H_\eps$, we have
\begin{multline}\label{eq:d2f}
D^2 f_\eps(u,v)[(h,h'),(k,k')]=
D^2 f_\eps^J (u)[h,k]
+D^2 f_\eps^K (v)[h',k']
\\
-\b\int_{\Omega_\varepsilon} v^2 h k
-2 \b\int_{\Omega_\varepsilon} u v h k'
-2 \b\int_{\Omega_\varepsilon} u v h' k
-\b\int_{\Omega_\varepsilon} u^2 h' k'.
\end{multline}

By \eqref{eq:UQ1}, if we set $a(Q)=\sqrt{J_1(Q)/J_2(Q)}$ and $b(Q)=\sqrt{J_1(Q)}$,
we have that $U^Q(x)=a(Q) W(b(Q) x)$ and so $U_P(x)=\chi(x-P) a(\eps P) W(b(\eps P)(x-P))$.
Therefore, we have:
\begin{align*}
\de_{P_i} U_P(x)
=&\;\de_{P_i} \left(\chi(x-P)  U^Q(x-P)\right)
\\
=& -  U^Q(x-P) \de_{x_i} \chi(x-P)
+\chi(x-P)  \de_{P_i} U^Q(x-P)
\\
=& -  U^Q(x-P) \de_{x_i} \chi(x-P)
+\eps \chi(x-P)\de_{P_i} a(\eps P) W(b(\eps P)(x-P))
\\
&+ \eps \chi(x-P)a(\eps P) \de_{P_i} a(\eps P) \n W(b(\eps P)(x-P))\cdot (x-P)
\\
&-\chi(x-P) a(\eps P) b(\eps P) (\de_{x_i} W)(b(\eps P)(x-P)).
\end{align*}
Hence
\begin{equation} \label{eq:de_iU}
\de_{P_i} U_P (x)=-\de_{x_i} U_P(x)+O(\eps).
\end{equation}
Analogously, we can prove that
\begin{equation} \label{eq:de_iV}
\de_{P_i} V_{P'} (x)=\de_{P'_i} V_{P'} (x)=-\de_{x_i} V_{P'}(x)+O(\eps).
\end{equation}
We recall that
\[
T_{(U_P,V_{P'})} Z^\eps = {\rm span}_{\H_\eps}
\{(\de_{P_1} U_P, \de_{P_1} V_{P'}) ,(\de_{P_2} U_P, \de_{P_2} V_{P'}) ,
(\de_{P_3}U_P, \de_{P_3} V_{P'}) \}.
\]
We set
\[
{\mathcal V}_\eps =
{\rm span}_{\H_\eps} \{(U_P,V_{P'}),(\de_{x_1} U_P, \de_{x_1} V_{P'}) ,
(\de_{x_2} U_P, \de_{x_2} V_{P'}) ,(\de_{x_3}U_P, \de_{x_3} V_{P'}) \}.
\]
By \eqref{eq:de_iU} and \eqref{eq:de_iV}, therefore it suffices to prove equation \eqref{eq:inv} for all
$(h,h')\in$ $ {\rm span}_{\H_\eps}\{(U_P,V_{P'}),(\phi, \phi')\}$, where $(\phi, \phi')$ is
orthogonal to ${\mathcal V}_\eps$. Precisely we shall prove that there exist $C_{1},C_{2}>0$ such that,
for all $\eps>0$ small enough, one has:
\begin{align}
(L_{\eps,Q}(U_P,V_{P'}) \mid (U_P,V_{P'})) \le & \; - C_{1}< 0,      \label{eq:neg}
\\
(L_{\eps,Q}(\phi,\phi') \mid (\phi,\phi')) \ge & \; C_{2} \|(\phi,\phi')\|^2, 
\qquad \textrm{for all } (\phi,\phi') \perp {\mathcal V}_\eps.        \label{eq:claim}
\end{align}

\vspace{0.2cm}
\noindent\textsc{Proof of \eqref{eq:neg}}. \quad By \eqref{eq:d2f}, we get:
\begin{multline}\label{eq:neg0}
D^2 f_\eps(U_P,V_{P'})[(U_P,V_{P'}),(U_P,V_{P'})]
\\
= D^2 f_\eps^J (U_P)[U_P,U_P]
+D^2 f_\eps^K (V_{P'})[V_{P'},V_{P'}]
-6\b\int_{\Omega_\varepsilon} U_P^2 V_{P'}^2.
\end{multline}
Let us study the first term of the right hand side of \eqref{eq:neg0}.
\begin{align*}
D^2 f_\eps^J (U_P)[U_P,U_P]
=&\int_{\Omega_\varepsilon} |\nabla U_P|^2 
+\int_{\Omega_\varepsilon} J_1 (\eps x) U_P^2 
-3\int_{\Omega_\varepsilon} J_2 (\eps x) U_P^{4} 
\\
=& \int\limits_{(\Omega-Q)/\eps\,\cap B_{\eps^{-1/4}}} \!\!\!\!\!\!\!
\left[ |\nabla U^Q|^2
+ J_1 (\eps x+Q) \left(U^Q\right) ^2 \!\!
-3 J_2 (\eps x+Q) \left(U^Q\right)^{4} \right]
\\
&+o(\eps)
\\
=&\int_{\RT}
\left[|\nabla U^Q|^2
+ J_1(Q) \left(U^Q\right) ^2
-3 J_2(Q) \left(U^Q\right)^{4} \right]
\\
&+\int_{\RT} \big( J_1(\eps x+Q) -J_1(Q)\big) \left(U^Q\right) ^2 
\\
&-3 \int_{\RT} \big( J_2(\eps x+Q) -J_2(Q)\big) \left(U^Q\right) ^4 
+o(\eps)
\\
=&\; -2 \int_{\RT} J_2(Q) \left(U^Q\right)^{4}  +O(\eps)
\\
=&\; -2 J_1(Q)^{\frac 12}J_2(Q)^{-1} \int_{\RT} W^{4}  +O(\eps)
\le -c_1.
\end{align*}
In a similar way it is possible to prove that
\[
D^2 f_\eps^K (V_{P'})[V_{P'},V_{P'}]\le -c_2.
\]
Finally, by Lemma \ref{le:UV}, we know that
\[
\int_{\Omega_\varepsilon} U_P^2 V_{P'}^2=o(\eps),
\]
and so equation \eqref{eq:neg} is proved.

\vspace{0.2cm}
\noindent\textsc{Proof of \eqref{eq:claim}}. \quad
Recalling the definition of $\chi$, (see \eqref{eq:chi}), we set $\chi_1 := \chi$ and
$ \chi_{2}:= 1-\chi_{1}$.
Given $(\phi, \phi') \perp {\mathcal V}_\eps$, let us consider the functions
\begin{align}
&\phi_{i}(x)=\chi_{i}(x-P)\phi(x),\quad i=1,2; \label{eq:phichi}
\\
&\phi_{i}'(x)=\chi_{i}(x-P')\phi'(x),\quad i=1,2. \label{eq:phichi'}
\end{align}
With calculations similar to those of \cite{AMS}, we have
\begin{align}
\| \phi \|^2 =&\; \| \phi_1 \|^2 + \| \phi_2 \|^2 +
\underbrace{2\int_{\Omega_\varepsilon} \chi_{1}\chi_{2}(\phi^{2}+|\nabla \phi|^{2})}_{I_\phi}
+ O(\eps^{1/4})\| \phi \|^2 \label{eq:phi},
\\
\| \phi' \|^2 =&\; \| \phi'_1 \|^2 + \| \phi'_2 \|^2 +
\underbrace{2\int_{\Omega_\varepsilon} \chi_{1}\chi_{2}((\phi')^{2}+|\nabla \phi'|^{2})}_{I_{\phi'}}
+ O(\eps^{1/4})\| \phi' \|^2 \label{eq:phi'}.
\end{align}
We need to evaluate the three terms in the equation below:
\begin{align}
(L_{\eps,Q}(\phi,\phi') \mid (\phi,\phi'))=&\;
(L_{\eps,Q}(\phi_{1},\phi'_1) \mid (\phi_{1},\phi'_1))+
(L_{\eps,Q}(\phi_{2},\phi'_2) \mid (\phi_{2},\phi'_2)) \nonumber
\\
&+ 2(L_{\eps,Q}(\phi_{1},\phi'_1) \mid (\phi_{2},\phi'_2)).\label{eq:L}
\end{align}

Let us start with $(L_{\eps,Q}(\phi_{1},\phi'_1) \mid (\phi_{1},\phi'_1))$.
Since $\b<0$, we get
\begin{align}
(L_{\eps,Q}(\phi_{1},\phi'_1) \mid (\phi_{1},\phi'_1)) =&\;
D^2 f_\eps^J (U_P)[\phi_1,\phi_1]
+D^2 f_\eps^K (V_{P'})[\phi'_1,\phi'_1] \nonumber
\\
& -4\b\!\int_{\Omega_\varepsilon}\!\!\! U_P V_{P'} \phi_1 \phi'_1
-\b\!\int_{\Omega_\varepsilon}\!\!\! U_P^2 {\phi'_1}^2
-\b\!\int_{\Omega_\varepsilon}\!\!\! V_{P'}^2 \phi_1^2 \nonumber
\\
>&\; D^2 f_\eps^J (U_P)[\phi_1,\phi_1]
+D^2 f_\eps^K (V_{P'})[\phi'_1,\phi'_1] \nonumber
\\
& -4\b\!\int_{\Omega_\varepsilon}\!\!\! U_P V_{P'} \phi_1 \phi'_1.        \label{eq:L11}
\end{align}
Arguing as in Lemma \ref{le:UV}, we know that
\begin{equation}\label{eq:L11-UV}
\int_{\Omega_\varepsilon} U_P V_{P'} \phi_1 \phi'_1=o(\eps).
\end{equation}
Therefore we need only to study the first two terms of the right hand side of
\eqref{eq:L11}. For simplicity, we can assume that $Q=\eps P$ is the origin
$\mathcal O$. In this case, we recall that we denote with $U^{\mathcal O}$ the
unique solution of \eqref{eq:UQ} whenever $Q={\mathcal O}$, while we denote
with $U_{\mathcal O}$ the truncation of $U^{\mathcal O}$, namely 
$U_{\mathcal O}=\chi \, U^{\mathcal O}$, where $\chi$ is defined in \eqref{eq:chi}. 
We have
\begin{align*}
D^2 f_\eps^J (U_{\mathcal O})[\phi_1,\phi_1]
=&\;\int_{\Omega_\varepsilon}
\left[ |\nabla \phi_1|^2
+ J_1(\eps x)\phi_1^2
-3 J_2(\eps x) U_{\mathcal O}^{2} \phi_1^2
\right]
\\
=&\;\int_{\RT}
\left[ |\nabla \phi_1|^2
+ J_1(\eps x)\phi_1^2
- 3 J_2(\eps x)\left(U^{\mathcal O}\right)^{2} \phi_1^2
\right] 
+o(\eps) \|\phi\|^2
\\
=&\; D^2 F^{J(\mathcal O)} (U^{\mathcal O})[\phi_1,\phi_1]
\\
&+\int_{\RT} \big( J_1(\eps x) -J_1(\mathcal O)\big) \phi_1^2 
\\
&- 3 \int_{\RT} \big( J_2(\eps x) -J_2(\mathcal O)\big) \left(U^{\mathcal O}\right)^{2}\phi_1^2 
+o(\eps) \|\phi\|^2
\\
\ge &\; D^2 F^{J(\mathcal O)} (U^{\mathcal O})[\phi_1,\phi_1]
- c_3 \,\eps \int_{\RT} |x| \, \phi_1^2 
+O(\eps) \|\phi\|^2
\\
=&\; D^2 F^{J(\mathcal O)} (U^{\mathcal O})[\phi_1,\phi_1]
+O(\eps^{3/4}) \|\phi\|^2,
\end{align*}
therefore
\begin{equation}\label{eq:F11}
D^2 f_\eps^J (U_{\mathcal O})[\phi_1,\phi_1]
\ge D^2 F^{J(\mathcal O)} (U^{\mathcal O})[\phi_1,\phi_1]
+O(\eps^{3/4}) \|\phi\|^2.
\end{equation}
We recall that $\phi$ is orthogonal to
\[
{{\mathcal V}^U_\eps} =
{\rm span}_{H^1_0(\Omega_\varepsilon)}
\{U_{\mathcal O}, \de_{x_1} U_{\mathcal O},\de_{x_2} U_{\mathcal O},\de_{x_3} U_{\mathcal O} \}.
\]
Moreover by \cite{C}, we know that if $\tf$ is orthogonal to $\mathcal V$ with
\[
{\mathcal V}^U =
{\rm span}_{H^1(\RT)} \{U^{\mathcal O}, \de_{x_1} U^{\mathcal O}, 
\de_{x_2} U^{\mathcal O}, \de_{x_3}U^{\mathcal O} \},
\]
then the fact that $U^{\mathcal O}$ is a Mountain Pass critical point 
of $F^{J(\mathcal O)}$ implies that
\begin{equation}\label{eq:D2F+}
D^2 F^{J(\mathcal O)}(U^{\mathcal O})[\tilde \phi,\tilde \phi]>c_4 \|\tilde \phi\|^2_{\RT}
\qquad \textrm{for all } \tilde \phi \perp {\mathcal V}^U.
\end{equation}
We can write $\phi_1=\xi + \zeta$, where $\xi \in {\mathcal V}^U$ and 
$\zeta \perp {\mathcal V}^U$. More precisely
\[
\xi=(\phi_1  \mid  U^{\mathcal O})_{\RT} \, U^{\mathcal O} \|U^{\mathcal O}\|^{-2}_{\RT}
+\sum_{i=1}^{3}(\phi_1  \mid  \de_{x_i} U^{\mathcal O})_{\RT} \,
\de_{x_i} U^{\mathcal O} \|\de_{x_i} U^{\mathcal O}\|^{-2}_{\RT}.
\]
Let us calculate $(\phi_1 \mid U^{\mathcal O})_{\RT}$. By the
exponential decay of $U^{\mathcal O}$ and since $\phi \perp
{\mathcal V}_\eps^U$, we have
\begin{align*}
(\phi_1  \mid  U^{\mathcal O})_{\RT}
=&\; \int_{\RT} \nabla \phi_1 \cdot \nabla  U^{\mathcal O}
+\int_{\RT} \phi_1  U^{\mathcal O}
\\
=&\;\int_{\Omega_\varepsilon} \nabla \phi_1 \cdot \nabla  U_{\mathcal O}
+\int_{\Omega_\varepsilon} \phi_1  U_{\mathcal O}
+o(\eps)\|\phi \|
\\
=&\;\int_{\Omega_\varepsilon} \nabla \phi \cdot \nabla  U_{\mathcal O}
+\int_{\Omega_\varepsilon} \phi \, U_{\mathcal O}
+o(\eps)\|\phi \|
=o(\eps)\|\phi \|.
\end{align*}
In a similar way, we can prove also that 
$(\phi_1  \mid \de_{x_i} U^{\mathcal O})_{\RT}=o(\eps)\|\phi \|$,
and so
\begin{align}
&\|\xi\|_{\RT} =o(\eps)\|\phi \|,    \label{eq:xi}
\\
&\|\zeta\|_{\RT} =\|\phi_1 \| + o(\eps)\|\phi \|.     \label{eq:zeta}
\end{align}
Let us estimate $D^2 F^{J(\mathcal O)}(U^{\mathcal O})[\phi_1,\phi_1]$. We get:
\begin{align}
D^2 F^{J(\mathcal O)}(U^{\mathcal O})[\phi_1,\phi_1]
=&\;D^2 F^{J(\mathcal O)}(U^{\mathcal O})[\zeta,\zeta]
+2 D^2 F^{J(\mathcal O)}(U^{\mathcal O})[\zeta,\xi] \nonumber
\\
&+D^2 F^{J(\mathcal O)}(U^{\mathcal O})[\xi,\xi].\label{eq:D2F}
\end{align}
By \eqref{eq:D2F+} and \eqref{eq:zeta}, since $\zeta \perp {\mathcal V}^U$, we know that
\[
D^2 F^{J(\mathcal O)}(U^{\mathcal O})[\zeta,\zeta]>c_3 \|\zeta\|^2_{\RT}
=c_3 \|\phi_1 \|^2 + o(\eps)\|\phi \|^2,
\]
while, by \eqref{eq:xi} and straightforward calculations, we have
\begin{align*}
&D^2 F^{J(\mathcal O)}(U^{\mathcal O})[\zeta,\xi] = o(\eps)\|\phi \|^2,
\\
&D^2 F^{J(\mathcal O)}(U^{\mathcal O})[\xi,\xi]  = o(\eps)\|\phi \|^2.
\end{align*}
By these last two estimates, \eqref{eq:D2F} and \eqref{eq:F11}, we can say that
\[
D^2 f_\eps^J (U_{\mathcal O})[\phi_1,\phi_1]
>c_4 \|\phi_1 \|^2 + O(\eps^{3/4})\|\phi\|^2.
\]
Hence, in the general case, we infer that, for all $Q\in \Omega_0$ with $|Q|\le\mu$,
\begin{equation}\label{eq:F11-Ufin}
D^2 f_\eps^J (U_P)[\phi_1,\phi_1]
>c_4 \|\phi_1 \|^2 + O(\eps^{3/4})\|\phi\|^2,
\end{equation}
and, analogously,
\begin{equation}\label{eq:F11-Vfin}
D^2 f_\eps^K (V_{P'})[\phi'_1,\phi'_1]
>c_5 \|\phi'_1 \|^2 + O(\eps^{1/2})\|\phi'\|^2.
\end{equation}
By \eqref{eq:L11}, \eqref{eq:L11-UV}, \eqref{eq:F11-Ufin} and \eqref{eq:F11-Vfin}, we can say that
\begin{equation}\label{eq:L11-fin}
(L_{\eps,Q}(\phi_{1}, \phi'_1) \mid (\phi_{1},\phi'_1))>
c_6 \|(\phi_1,\phi'_1) \|^2
+ O(\eps^{1/2})\|(\phi,\phi'_1)\|^2.
\end{equation}

Let us now evaluate $(L_{\eps,Q}(\phi_{2},\phi'_2) \mid (\phi_{2},\phi'_2))$.
Arguing as in Lemma \ref{le:UV}, since $\b<0$
and using the definition of $\chi_i$ and the exponential decay of
$U_P$ and of $V_{P'}$, we easily get:
\begin{align}
(L_{\eps,Q}(\phi_{2},\phi'_2) \mid (\phi_{2},\phi'_2))
=&\;
D^2 f_\eps^J (U_P)[\phi_2,\phi_2]
+D^2 f_\eps^K (V_{P'})[\phi'_2,\phi'_2] \nonumber
\\
& -4 \b \int_{\Omega_\varepsilon} U_P V_{P'} \phi_2 \phi'_2
-\b \int_{\Omega_\varepsilon} U_P^2 {\phi'_2}^2
-\b \int_{\Omega_\varepsilon} V_{P'}^2 \phi_2^2 \nonumber
\\
\ge &\;
D^2 f_\eps^J (U_P)[\phi_2,\phi_2]
+D^2 f_\eps^K (V_{P'})[\phi'_2,\phi'_2]
+o(\eps)\|(\phi,\phi')\|^{2} \nonumber
\\
\ge & \;c_7 \|(\phi_{2},\phi'_2)\|^{2}+o(\eps)\|(\phi,\phi')\|^{2}. \label{eq:L22}
\end{align}

Let us now study $(L_{\eps,Q}(\phi_{1},\phi'_1) \mid (\phi_{2},\phi'_2))$.
Arguing as in Lemma \ref{le:UV}, we get
\begin{align}
(L_{\eps,Q}(\phi_{1},\phi'_1) \mid (\phi_{2},\phi'_2))
=&\;
D^2 f_\eps^J (U_P)[\phi_1,\phi_2]
+D^2 f_\eps^K (V_{P'})[\phi'_1,\phi'_2] \nonumber
\\
& -2 \b \int_{\Omega_\varepsilon} U_P V_{P'} \phi_1 \phi'_2
-2 \b \int_{\Omega_\varepsilon} U_P V_{P'} \phi_2 \phi'_1 \nonumber
\\
& - \b \int_{\Omega_\varepsilon} U_P^2 \phi'_1 \phi'_2
- \b \int_{\Omega_\varepsilon} V_{P'}^2 \phi_1 \phi_2 \nonumber
\\
= & \;
D^2 f_\eps^J (U_P)[\phi_1,\phi_2]
+D^2 f_\eps^K (V_{P'})[\phi'_1,\phi'_2] \nonumber
\\
& -\b \int_{\Omega_\varepsilon} U_P^2 \phi'_1 \phi'_2
-\b \int_{\Omega_\varepsilon} V_{P'}^2 \phi_1 \phi_2 +o(\eps)\|(\phi,\phi')\|^{2}. \label{eq:L12}
\end{align}
Using the definition of $\chi_i$ and the exponential decay of
$U_P$ and of $V_{P'}$, we easily get:
\begin{align}
&D^2 f_\eps^J (U_P)[\phi_1,\phi_2] \ge
c_8 I_\phi + O(\eps^{1/4})\|\phi\|^{2}, \label{eq:L12-fU}
\\
& D^2 f_\eps^K (V_{P'})[\phi'_1,\phi'_2] \ge
c_9 I_{\phi'} + O(\eps^{1/4})\|\phi'\|^{2}, \label{eq:L12-fV}
\end{align}
where $I_\phi$ and $I_{\phi'}$ are defined, respectively in \eqref{eq:phi} and
\eqref{eq:phi'}.
Moreover, by the definition of $\chi$, (see \eqref{eq:chi}), and by the 
definitions of $\phi_i$ and $\phi_i'$, (see \eqref{eq:phichi} and \eqref{eq:phichi'}),
\[
\phi_1 (x) \,\phi_2(x)
=\chi(x-P)(1-\chi(x-P))\phi^2(x)\ge 0, \qquad \textrm{for all } x\in \RT,
\]
and so, also
\[
\phi'_1 (x) \,\phi'_2(x)\ge 0,\qquad \textrm{for all } x\in \RT.
\]
Therefore
\begin{equation}\label{eq:L12-U-V}
-\b\int_{\Omega_\varepsilon} U_P^2 \phi'_1 \phi'_2
-\b\int_{\Omega_\varepsilon} V_{P'}^2 \phi_1 \phi_2
\ge 0.
\end{equation}
By \eqref{eq:L12}, \eqref{eq:L12-fU}, \eqref{eq:L12-fV} and \eqref{eq:L12-U-V}, we infer
\begin{equation}\label{eq:L12-fin}
(L_{\eps,Q}(\phi_{1},\phi'_1) \mid (\phi_{2},\phi'_2))
\ge c_{10} \left(I_{\phi}+ I_{\phi'}\right)
+ O(\eps^{1/4})\|(\phi, \phi')\|^{2}.
\end{equation}
Hence, by \eqref{eq:L}, \eqref{eq:L11-fin}, \eqref{eq:L22}, \eqref{eq:L12-fin} and
recalling \eqref{eq:phi} and \eqref{eq:phi'}, we get
\[
(L_{\eps,Q}(\phi,\phi') \mid (\phi,\phi'))
\ge c_{11} \|(\phi,\phi')\|^{2}+ O(\eps^{1/4})\|(\phi,\phi')\|^{2}.
\]
This completes the proof of the lemma.
\end{proof}

\section{The finite dimensional reduction}

By means of the Liapunov-Schmidt reduction, the existence of critical points of $f_\eps$ can be reduced
to the search of critical points of an auxiliary finite-dimensional
functional. 

\begin{lemma}\label{lem:w}
Fix $\mu>0$. For $\eps>0$ small enough and for all $Q\in \Omega_0$ with $|Q|\le \mu$, there exists a unique
$(w,w')=\big(w(\eps,Q),w'(\eps,Q)\big) \in \H_\eps$ of class $C^{1}$ such that:
\begin{enumerate}
\item $\big(w(\eps, Q),w'(\eps,Q)\big)\in (T_{(U_P,V_{P'})} Z^\eps)^{\perp}$;
\item $\nabla f_\eps (U_P + w,V_{P'} +w')\in T_{(U_P,V_{P'})} Z^\eps$.
\end{enumerate}
Moreover, the functional $\A_\eps \colon \Omega_0 \to \R$, defined as:
\[
\A_\eps (Q):=f_\eps \big(U_{Q/\eps} +w(\eps,Q),V_{(Q+e_1\sqrt{\eps})/\eps}+w'(\eps,Q)\big)
\]
is of class $C^1$ and satisfies:
\[
\nabla \A_\eps(Q_0)=0
\quad \Longleftrightarrow \quad
\nabla f_\eps\left(U_{Q_0/\eps}+w(\eps,Q_0),V_{(Q_0+e_1 \sqrt{\eps})/\eps}+ w'(\eps,Q_0)\right)=0.
\]
\end{lemma}

\begin{proof}
Let $\P=\P_{\eps, Q}$ denote the projection onto $(T_{(U_P,V_{P'})} Z^\eps)^\perp$. We want
to find a solution $(w,w') \in (T_{(U_P,V_{P'})} Z^\eps)^{\perp}$ of the equation
\[
P\nabla f_\eps(U_P +w,V_{P'}+w')=0.
\]
One has that
\[
\nabla f_\eps(U_P+w,V_{P'}+w')=
\nabla f_\eps (U_P,V_{P'})+D^2 f_\eps(U_P,V_{P'})[w,w']+R(U_P,V_{P'},w,w')
\]
with $\|R(U_P,V_{P'},w,w')\|=o(\|(w,w')\|)$, uniformly
with respect to $(U_P,V_{P'})$. Therefore, our equation is:
\begin{equation}\label{eq:eq-w}
L_{\eps,Q}(w,w') + \P\nabla f_\eps (U_P,V_{P'})+\P R(U_P,V_{P'},w,w')=0.
\end{equation}
According to Lemma \ref{le:inv}, this is equivalent to
\[
(w,w') = N_{\eps,Q}(w,w'),
\]
where
\[
N_{\eps,Q}(w,w')=-\left( L_{\eps,Q}\right)^{-1}\big( \P \nabla f_\eps (U_P,V_{P'})
+\P R(U_P,V_{P'},w,w')\big).
\]
By \eqref{eq:nf} it follows that
\begin{equation}\label{eq:N}
\|N_{\eps,Q}(w,w')\| = O(\eps^{1/2}) + o(\|(w,w')\|).
\end{equation}
Therefore it is easy to check that $N_{\eps,Q}$ is a contraction on some ball in
$(T_{(U_P,V_{P'})} Z^\eps)^{\perp}$
provided that $\eps>0$ is small enough.
Then there exists a unique $(w,w')$ such that $(w,w')=N_{\eps,Q}(w,w')$.  Let us
point out that we cannot use the Implicit Function Theorem to find
$\big(w(\eps,Q),w'(\eps,Q)\big)$, because the map
$(\eps,u,v)\mapsto \P\nabla f_\eps (u,v)$ fails to be
$C^2$.  However, fixed $\eps>0$ small, we can apply the Implicit
Function Theorem to the map $(Q,w,w')\mapsto \P \nabla f_\eps (U_P + w,V_{P'}+w')$.
Then, in particular, the function $\big(w(\eps,Q),w'(\eps,Q)\big)$ turns out to be of class
$C^1$ with respect to $Q$.  Finally, it is a standard argument, see
\cite{AB,ABC}, to check that the critical points of $\A_\eps(Q)=f_\eps (U_P+w,V_{P'}+w')$
give rise to critical points of $f_\eps$.
\end{proof}

\begin{remark}\label{rem:w}
From (\ref{eq:N}) it immediately follows that:
\begin{equation}\label{eq:w}
\|(w,w')\|=O(\eps^{1/2}).
\end{equation}
\end{remark}

Let us now make the asymptotic expansion of the finite dimensional functional.

\begin{theorem}\label{th:sviluppo}
Fix $\mu>0$ and let $Q\in \Omega_0$ with $|Q|\le \mu$, $Q'=Q+ \sqrt\eps \, e_1$, $P=Q/\eps \in \Omega_\varepsilon$ and 
$P'=Q'/\eps \in \Omega_\varepsilon$. 
Suppose {\bf (J)} and {\bf (K)}. Then,
for $\eps$ sufficiently small, we get:
\begin{equation}\label{eq:A}
\A_\eps (Q)=
f_\eps\big(U_P + w(\eps,Q),V_{P'} + w'(\eps,Q)\big)
= c_0 \G(Q) +o(\eps^{1/4}),
\end{equation}
where $\G \colon \Omega_0 \to \R$ is defined in \eqref{eq:Gamma}, namely
\[
\G(Q)=J_1(Q)^{\frac 12}J_2(Q)^{-1} + K_1(Q)^{\frac 12}K_2(Q)^{-1},
\]
and
\begin{equation}\label{eq:c0}
c_0 := \frac{1}{2}\int_{\RT} W^4
\end{equation}
with $W$ the unique solution of \eqref{eq:W}.
\end{theorem}

\begin{proof}
We have:
\begin{align*}
\A_\eps (Q)=&
f_\eps\big(U_P + w(\eps,Q),V_{P'} + w'(\eps,Q)\big)
\\
=& \; \frac{1}{2}\int_{\Omega_\varepsilon} |\nabla (U_P +w)|^2 +
\frac{1}{2}\int_{\Omega_\varepsilon} J_1(\eps x)(U_P+w)^2
-\frac{1}{4}\int_{\Omega_\varepsilon} J_2(\eps x)(U_P+w)^{4}
\\
& + \frac{1}{2}\!\int_{\Omega_\varepsilon} \!|\nabla (V_{P'}+w')|^2 +
\frac{1}{2}\!\int_{\Omega_\varepsilon} \!K_1(\eps x) (V_{P'}+w')^2
-\frac{1}{4}\!\int_{\Omega_\varepsilon} \!K_2(\eps x)(V_{P'}+w')^{4}
\\
& - \frac{\b}{2} \int_{\Omega_\varepsilon} (U_{P}+w)^2 (V_{P'}+w')^2.
\end{align*}
Therefore, by \eqref{eq:w} and Lemma \ref{le:UV},
\begin{align}
\A_\eps (Q)
=& \; \frac{1}{2}\int_{\Omega_\varepsilon} |\nabla U_P|^2 +
\frac{1}{2}\int_{\Omega_\varepsilon} J_1(\eps x) U_P^2
-\frac{1}{4}\int_{\Omega_\varepsilon} J_2(\eps x) U_P^{4} \nonumber
\\
&+ \frac{1}{2}\int_{\Omega_\varepsilon} |\nabla V_{P'}|^2 +
\frac{1}{2}\int_{\Omega_\varepsilon} K_1(\eps x) V_{P'}^2
-\frac{1}{4}\int_{\Omega_\varepsilon} K_2(\eps x) V_{P'}^{4}
+O(\eps^{1/2}) \nonumber
\\
=&\;
f_\eps^J (U_P)
+f_\eps^K (V_{P'})
+O(\eps^{1/2}).          \label{eq:A-tot}
\end{align}
Let us study the first term of the right hand side of \eqref{eq:A-tot}.
\begin{align*}
f_\eps^J(U_P)
=&\; \frac{1}{2}\int_{\Omega_\varepsilon} |\nabla U_P|^2 +
\frac{1}{2}\int_{\Omega_\varepsilon} J_1(\eps x) U_P^2
-\frac{1}{4}\int_{\Omega_\varepsilon} J_2(\eps x) U_P^{4}
\\
=&\; \frac 12 \int\limits_{(\Omega-Q)/\eps\,\cap B_{\eps^{-1/4}}} \!\!\!\!
\left[ |\nabla U^Q|^2 \!\!
+\! J_1(\eps x+Q) \left(U^Q\right) ^2 \!\!
-\frac 12 J_2(\eps x+Q)\left(U^Q\right)^{4} \right]\!
+o(\eps)
\end{align*}
\begin{align*}
=&\; \frac 12 \int_{\RT} |\nabla U^Q|^2
+\frac 12 \int_{\RT}  J_1(Q) \left(U^Q\right) ^2
-\frac 14 \int_{\RT}  J_2(Q)\left(U^Q\right)^{4}
\\
& + \frac 12 \int_{\RT} \big( J_1(\eps x+Q) -J_1(Q)\big) \left(U^Q\right) ^2 
\\
& - \frac 14 \int_{\RT} \big( J_2(\eps x+Q) -J_2(Q)\big) \left(U^Q\right) ^4 
+o(\eps)
\\
=&\; \frac 12  \int_{\RT}J_2(Q) \left(U^Q\right)^{4}  +o(\eps^{1/4})
\\
=&\; \frac 12  J_1(Q)^{\frac 12}J_2(Q)^{-1} \int_{\RT} W^{4}  +o(\eps^{1/4}).
\end{align*}
Hence
\begin{equation}\label{eq:A-U}
f_\eps^J(U_P)
=\frac 12  J_1(Q)^{\frac 12}J_2(Q)^{-1} \int_{\RT} W^{4}  +o(\eps^{1/4}).
\end{equation}
Analogously,
\begin{align*}
f_\eps^K(V_{P'})
=&\; \frac{1}{2}\int_{\Omega_\varepsilon} |\nabla V_{P'}|^2 +
\frac{1}{2}\int_{\Omega_\varepsilon} K_1(\eps x)V_{P'}^2
-\frac{1}{4}\int_{\Omega_\varepsilon} K_2(\eps x) V_{P'}^{4}
\\
=&\; \frac 12 \int\limits_{(\Omega-Q')/\eps\,\cap B_{\eps^{-1/4}}} \!\!
\left[ |\nabla V^Q|^2
+ K_1(\eps x+Q') \left(V^Q\right) ^2\right]
\\
& \; -\frac 14 \int\limits_{(\Omega-Q')/\eps\,\cap B_{\eps^{-1/4}}} \!\!
K_2(\eps x+Q')\left(V^Q\right)^{4}
+o(\eps)
\\
=&\; \frac 12 \int_{\RT} |\nabla V^Q|^2
+\frac 12 \int_{\RT}  K_1(Q) \left(V^Q\right) ^2
-\frac 14 \int_{\RT}  K_2(Q) \left(V^Q\right)^{4}
\\
& + \frac 12 \int_{\RT} \left( K_1(\eps x+Q+\sqrt{\eps} \, e_1) -K_1(Q)\right) \left(V^Q\right) ^2 
\\
& - \frac 14 \int_{\RT} \left( K_2(\eps x+Q+\sqrt{\eps} \, e_1) -K_2(Q)\right) \left(V^Q\right) ^4 
+o(\eps)
\\
=&\; \frac 12  \int_{\RT}  K_2(Q)\left(V^Q\right)^{4}  +o(\eps^{1/4})
\\
=&\; \frac 12  K_1(Q)^{\frac 12}K_2(Q)^{-1} \int_{\RT} W^{4}  +o(\eps^{1/4}).
\end{align*}
Therefore
\begin{equation}\label{eq:A-V}
f_\eps^J(V_{P'})
=\frac 12  K_1(Q)^{\frac 12}K_2(Q)^{-1} \int_{\RT} W^{4}  +o(\eps^{1/4}).
\end{equation}
Now \eqref{eq:A} follows immediately by \eqref{eq:A-tot}, \eqref{eq:A-U} and \eqref{eq:A-V}.
\end{proof}

\section{A multiplicity result and proofs of theorems}

In this section we give the proofs of our theorems. First of all, 
let us prove Theorem \ref{th1} as an easy 
consequence of the following multiplicity result:
\begin{theorem}\label{th:cat}
Let {\bf (J)} and {\bf (K)} hold and suppose $\G$ has
a compact set $X \subset \Omega_0$ where $\G$ achieves a strict local minimum (resp. maximum),
in the sense that there exist $\delta>0$ and a $\d$-neighborhood $X_\d \subset \Omega_0$ of $X$ such that
\[
b:=\inf\{\G(Q) : Q\in \partial X_{\d}\}> a:= \G_{|_X}, \quad
\left({\rm resp. }\; \sup\{\G(Q) : Q\in \partial X_{\d}\}< \G_{|_X} \right).
\]

Then there exists $\bar \eps>0$ such that \eqref{eq}
has at least $\cat(X,X_\d)$ solutions that concentrate near points of $X_{\d}$, provided
$\eps \in (0,\bar \eps)$. Here $\cat(X,X_\d)$ denotes the Lusternik-Schnirelman category of $X$ with
respect to $X_\d$.
\end{theorem}

\begin{proof}
First of all, we fix $\mu>0$ in such a way that $|Q|<\mu$ for all $Q\in X$.  We will apply the
finite dimensional procedure with such $\mu$ fixed. 

We will treat only the case of minima, being the other one similar.
We set $Y=\{Q\in X_{\d} :\A_{\eps}(Q)\le c_0 (a+b)/2\}$, being $c_0$ defined
in \eqref{eq:c0}. 
By \eqref{eq:A} it follows that there exists $\bar \eps>0$ such that
\begin{equation}\label{eq:X}
X  \subset Y \subset X_{\d},
\end{equation}
provided $\eps\in (0,\bar \eps)$. Moreover, if $Q \in \partial X_{\d}$ then
$\G(Q)\ge b$ and hence
\[
\A_{\eps}(Q)\ge c_0 \G (Q) + o(\eps^{1/4}) \ge  c_0 b + o(\eps^{1/4}).
\]
On the other side, if $Q\in Y$
then $\A_{\eps}(Q)\le c_0 (a+ b)/2$.  Hence, for $\eps$ small,
$Y$ cannot meet $\partial X_{\d}$
and this readily implies that $Y$ is compact.
Then $\A_{\eps}$ possesses at least $\cat(Y,X_{\d})$
critical points in $ X_{\d}$.  Using (\ref{eq:X}) and the properties of
the category one gets
\[
\cat(Y,Y)\ge \cat(X,X_{\d}).
\]
Moreover, by Lemma \ref{lem:w}, we know that to critical points of $\A_\eps$ 
there correspond critical points of $f_\eps$ and so solutions of \eqref{eqe}.
Let $Q_\eps \in X$ be one of these critical points, if $Q'_\eps= Q_\eps+\sqrt{\eps}\, e_1$, then
\[
(u_\eps^{Q_\eps},v_\eps^{Q_\eps})
=\left(U_{Q_\eps/\eps}+w(\eps, Q_\eps),V_{Q'_\eps/\eps}+w'(\eps, Q_\eps)\right)
\]
is a solution of \eqref{eqe}. Therefore
\begin{align*}
&u_\eps^{Q_\eps}(x/\eps)\simeq U_{Q_\eps/\eps}(x/\eps)= U^{Q_\eps}\left(\frac{x-Q_\eps}{\eps} \right)
\\
&v_\eps^{Q_\eps}(x/\eps)\simeq V_{Q'_\eps/\eps}(x/\eps)= V^{Q_\eps}\left(\frac{x-Q'_\eps}{\eps} \right)
\end{align*}
is a solution of \eqref{eq} and also the concentration result follows.
\end{proof}

Let us now give a short proof of Theorem \ref{thN}.

\begin{proofN}
We need only to observe that, in this case, the solutions of \eqref{eqN} 
will be found near $(\bar U^Q, \bar V^Q)$, properly
truncated, where $\bar U^Q$ is the unique solution of
\begin{equation*}
\left\{
\begin{array}{ll}
-\Delta  u+ J_1(Q) u= J_2(Q)u^{2p-1} &   \text{in }\RN,
\\
u >0 &   \text{in }\RN,
\\
u(0)=\max_{\RN} u,
\end{array}
\right.
\end{equation*}
and $\bar V_{Q}$ is the unique solution of
\begin{equation*}
\left\{
\begin{array}{ll}
-\Delta  v+ K_1(Q) v= K_2(Q)v^{2p-1} &   \text{in }\RN,
\\
v >0 &   \text{in }\R^{N},
\\
v(0)=\max_{\R^{N}} v,
\end{array}
\right.
\end{equation*}
for an appropriate choice of $Q \in \bar\Omega_0$. It is easy to see that
\begin{align*}
\bar U^Q(x)&=\big(J_1(Q)/J_2(Q)\big)^{1/(2p-2)} \cdot \bar W\left(\sqrt{J_1(Q)} \cdot x   \right),
\\
\bar V^{Q}(x)&=\big(K_1(Q)/K_2(Q)\big)^{1/(2p-2)} \cdot \bar W\left(\sqrt{K_1(Q)} \cdot x \right),
\end{align*}
where $\bar W$ is the unique solution of
\begin{equation*}
\left\{
\begin{array}
[c]{lll}
-\Delta  z+ z= z^{2p-1} &   \text{in }\R^{N},
\\
z>0 &   \text{in }\R^{N},
\\
z(0)=\max_{\R^{N}}z.
\end{array}
\right.
\end{equation*}
At this point, we can repeat the previous arguments, with suitable modifications.
\end{proofN}

\begin{remark}
Of course, the analogous of Theorem \ref{th:cat} holds also for problem \eqref{eqN}. 
\end{remark}

\end{document}